# Inifnite hypercomplex number system factorization methods

**Yakiv O. Kalinovsky,** Dr.Sc., Senior Researcher,
Institute for Information Recording National Academy of Science of Ukraine, Kyiv, Shpaka str. 2, 03113, Ukraine, E-mail: kalinovsky@i.ua

**Dmitry V. Lande,** Dr.Sc., Head of Department,
Institute for Information Recording National Academy of Science of Ukraine, Kyiv, Shpaka str. 2, 03113, Ukraine, E-mail: dwlande@gmail.com

**Yuliya E. Boyarinova,** PhD, Associate Professor,
National Technical University of Ukraine "KPI", Kyiv, Peremogy av. 37, 03056, Ukraine, E-mail: ub@ua.fm

**Iana V. Khitsko,** Junior Researcher,
National Technical University of Ukraine "KPI", Kyiv, Peremogy av. 37, 03056, Ukraine, E-mail: yannuary@yandex.ua

**Abstract**

The method of obtaining the set of noncanonical hypercomplex number systems by conversion of infinite hypercomplex number system to finite hypercomplex number system depending on multiplication rules and factorization method is described. Systems obtained by this method starting from the 3$^{rd}$ dimension are noncanonical. The obtained systems of even dimension can be re-factorized. As a result of it hypercomplex number system of two times less dimension are got.

## Introduction

Hypercomplex data presentation form is widely used for increasing algorythms performance in different areas of science and technology. Canonical hypercomplex number systems (HNS) are used more often; here belong complex, double and dual numbers, quaternions etc**[1]**. It is possible that different noncanonical hypercomplex number systems can be used either for some tasks for better algorythms performance or, as in cryptographic algorithms, to increase their complexity to prevent data hacking, depending on selected number systems **[2-3]**.

A number system is considered canonical if each multiplication of every basic element pair is equal to one of the basic elements with the coefficient from $\{-1;0;+1\}$. Providing that one of such basic elements multiplication is a sum of two or more summands and/or with the coefficient out of $\{-1;0;+1\}$, the HNS is called noncanonical**[1]**.

It is very important for noncanonical HNS applicability investigation to build the upper mentioned systems. In general, noncanonical HNS are built via iterations through structure constants of the basic elements in the multiplication table cells. Considering this, the constraint is set on the range values of structure constants, the unit element presence and type, and linear independence of basic elements. In addition, there is a method for getting the noncanonical HNS set that are isomorphic for the given number system via iterations through linear equation coefficients **[4-5]**.

## Factorization method of infinite hypercomplex number system

Let us consider the way of getting noncanonical hypercomplex number systems set, such as conversion the infinite HNS to finite HNS depending on specific multiplication rules and factorization method.

Let $\Gamma$ be discrete countable set with infinite basis $\{e_i\}, i = 1,....,n,....$ The involution operation* introduced for this set is as follows:

$$*: \Gamma \to \Gamma, \text{ where for } \forall e_i \in \Gamma \ e_i^* \in \Gamma. \qquad (1)$$

There is a basic element $e_1 \in \Gamma$, such that

$$e_1 \cdot e_i = e_i \cdot e_1 = e_i, \ (e_i^*)^* = e_i \qquad (2)$$

Multiplication operation (convolution) in the set is made according to the rule below:

$$e_i \cdot e_j = \sum_{k=1}^{\infty} C_{ij}^k e_k \qquad (3)$$

$\Gamma$ set is a hypercomplex system which can be both discrete and continuous under the conditions (1)-(3).

Conversion of the infinite HNS to finite HNS will be done considering the conditions of commutativity and the positivity of structure constants:

$$C_{ij}^k \geq 0; \quad C_{ii^*}^1 > 0; \quad C_{ij}^k = C_{ki^*}^j \geq 0. \qquad (4)$$

We will build the infinite hypercomplex system, from which subsequently factoring out the specific subgroup we will receive hypercomplex number systems of finite dimension.

It is given group $Z = \{-\infty, \infty\}$, on which a subgroup of automorphisms is selected $V = \{-1, 1\} \in Aut\ Z$, that for $\forall n \in Z$:

$$1(n) = n \in Z, -1(n) = -n \in Z. \qquad (5)$$

Next step is to factor $Z$ group by the subgroup $V$ and obtain a set $Z/V = \Gamma = N \cup \{0\}$. For this set the multiplication rules (convolution) shall be defined.

If convolution for $Z$ group is

$$\sigma_n \cdot \sigma_m = \sigma_{n+m}. \qquad (6)$$

or, taking into account $Z$ group features, we have the relation:

$$n \cdot m = n + m. \qquad (7)$$

Then convolution for $Z/V = \Gamma$ will be:

$$n \cdot m = \{n, -n\} \cdot \{m, -m\} = (n+m) + (m-n) + (-n-m) + (n-m) = \\ = \{(n+m), -(n+m)\} + \{(n-m), -(n-m)\} \qquad (8)$$

The relation holds

$$\frac{1}{2}(\sigma_n \cdot \sigma_{-n}) = \sigma_n. \qquad (9)$$

From the relations listed above:

$$\sigma_n \cdot \sigma_m = \frac{1}{2}(\sigma_n + \sigma_{-n}) \cdot \frac{1}{2}(\sigma_m + \sigma_{-m}) = \frac{1}{4}(\sigma_{n+m} + \sigma_{-(n+m)}) + (\sigma_{|n-m|} + \sigma_{-|n-m|}). \qquad (10)$$

Then the convolution (9) considering (10) will be:

$$\sigma_n \cdot \sigma_m = \frac{1}{2}(\sigma_{n+m} + \sigma_{|n-m|}). \qquad (11)$$

With the infinite hypercomplex system and convolution we select subgroups for which the factorization will be done.

Subgroups such as $\{A_2, A_3, A_4,...\}$ can be taken, where $A_2 = \{a_k : a_k = 2k, k \in \Gamma\}$, $A_3 = \{a_k : a_k = 3k, k \in \Gamma\}$, $A_4 = \{a_k : a_k = 4k, k \in \Gamma\}$ etc. It should be noted, that each subgroup is infinite, and the number of subgroups is countable.

**Analysis of obtained hypercomplex number systems**

As a result of factorization of the HNS $\Gamma = N \cup \{0\}$ by the subgroup $A_2 = \{a_k : a_k = 2k, k \in \Gamma\}$ with convolution (11), the finite hypercomplex number system of the second dimension $G_2 = \Gamma/A_2 = \{e_1, e_2\}$ is obtained, specifically, double number system. Multiplication table of this system can be seen below (Table 1):

Table 1.

| $e_1$ | $e_2$ |
|---|---|
| $e_2$ | $e_1$ |

Further, increasing the dimension, we obtain noncanonical HNS.

The result of factorized hypercomplex system $\Gamma = N \cup \{0\}$ by the subgroup $A_3 = \{a_k : a_k = 3k, k \in \Gamma\}$ with convolution (11), HNS $G_3 = \Gamma/A_3 = \{e_1, e_2, e_3\}$ of third dimension is obtained. Multiplication table of this system (Table 2):

Table 2.

| $e_1$ | $e_2$ | $e_3$ |
|---|---|---|
| $e_2$ | $\frac{1}{2}(e_1 + e_3)$ | $\frac{1}{2}(e_1 + e_2)$ |
| $e_3$ | $\frac{1}{2}(e_1 + e_2)$ | $\frac{1}{2}(e_1 + e_2)$ |

The result of factorized hypercomplex system $\Gamma = N \cup \{0\}$ by the subgroup $A_4 = \{a_k : a_k = 4k, k \in Q\}$ with convolution (11), HNS $G_4 = \Gamma/A_4 = \{e_1, e_2, e_3, e_4\}$ of forth dimension is built. Multiplication table of the obtained system is listed below (Table 3):

Table 3.

| $e_1$ | $e_2$ | $e_3$ | $e_4$ |
|---|---|---|---|
| $e_2$ | $\frac{1}{2}(e_1 + e_3)$ | $\frac{1}{2}(e_2 + e_4)$ | $\frac{1}{2}(e_1 + e_3)$ |
| $e_3$ | $\frac{1}{2}(e_2 + e_4)$ | $e_1$ | $e_2$ |
| $e_4$ | $\frac{1}{2}(e_1 + e_3)$ | $e_2$ | $\frac{1}{2}(e_1 + e_3)$ |

As far as there is a regularity in the multiplication table of HNS, listed above, such as $e_1 \cdot e_1 = e_1$, $e_3 \cdot e_3 = e_1$, we can do the re-factorization by the subset $\{e_1, e_3\}$. Re-factorization of $G_4/\{e_1, e_3\}$ led to the hypercomplex number system of second dimension with multiplication table (Table 4):

|  | $e_1$ | $e_2$ |
|---|---|---|
| $e_2$ |  | $\frac{1}{2}(e_1 + e_2)$ |

The result of factorized hypercomplex system $\Gamma = N \cup \{0\}$ by the subgroup $A_5 = \{a_k : a_k = 5k, k \in \Gamma\}$ with multiplication rules (11), HNS $G_5 = \Gamma/A_5 = \{e_1, e_2, e_3, e_4, e_5\}$ of dimension 5 is built. See the convolution of the obtained system in the table:

Table 5.

|  | $e_1$ | $e_2$ | $e_3$ | $e_4$ | $e_5$ |
|---|---|---|---|---|---|
| $e_2$ |  | $\frac{1}{2}(e_1 + e_3)$ | $\frac{1}{2}(e_2 + e_4)$ | $\frac{1}{2}(e_3 + e_5)$ | $\frac{1}{2}(e_1 + e_4)$ |
| $e_3$ |  | $\frac{1}{2}(e_2 + e_4)$ | $\frac{1}{2}(e_1 + e_5)$ | $\frac{1}{2}(e_1 + e_2)$ | $\frac{1}{2}(e_2 + e_3)$ |
| $e_4$ |  | $\frac{1}{2}(e_3 + e_5)$ | $\frac{1}{2}(e_1 + e_2)$ | $\frac{1}{2}(e_1 + e_2)$ | $\frac{1}{2}(e_2 + e_3)$ |
| $e_5$ |  | $\frac{1}{2}(e_1 + e_4)$ | $\frac{1}{2}(e_2 + e_3)$ | $\frac{1}{2}(e_2 + e_3)$ | $\frac{1}{2}(e_1 + e_4)$ |

Let us factor hypercomplex system $\Gamma = N \cup \{0\}$ by the subgroup $A_6 = \{a_k : a_k = 6k, k \in \Gamma\}$ the same way as by the previous subgroups and will get the HNS $G_6 = \Gamma/A_6 = \{e_1, e_2, e_3, e_4, e_5, e_6\}$ of dimension 6. The convolution of the obtained system will be as follows in Table 6:

Table 6.

|  | $e_1$ | $e_2$ | $e_3$ | $e_4$ | $e_5$ | $e_6$ |
|---|---|---|---|---|---|---|
| $e_2$ |  | $\frac{1}{2}(e_1 + e_3)$ | $\frac{1}{2}(e_2 + e_4)$ | $\frac{1}{2}(e_3 + e_5)$ | $\frac{1}{2}(e_4 + e_6)$ | $\frac{1}{2}(e_1 + e_5)$ |
| $e_3$ |  | $\frac{1}{2}(e_2 + e_4)$ | $\frac{1}{2}(e_1 + e_5)$ | $\frac{1}{2}(e_2 + e_6)$ | $\frac{1}{2}(e_1 + e_3)$ | $\frac{1}{2}(e_2 + e_4)$ |
| $e_4$ |  | $\frac{1}{2}(e_3 + e_5)$ | $\frac{1}{2}(e_2 + e_6)$ | $e_1$ | $e_2$ | $e_3$ |
| $e_5$ |  | $\frac{1}{2}(e_4 + e_6)$ | $\frac{1}{2}(e_1 + e_3)$ | $e_2$ | $\frac{1}{2}(e_1 + e_3)$ | $\frac{1}{2}(e_2 + e_4)$ |
| $e_6$ |  | $\frac{1}{2}(e_1 + e_5)$ | $\frac{1}{2}(e_2 + e_4)$ | $e_3$ | $\frac{1}{2}(e_2 + e_4)$ | $\frac{1}{2}(e_1 + e_5)$ |

Since there are regularities in multiplication the above HNS table, that $e_1 \cdot e_1 = e_4 \cdot e_4 = e_1$, $e_1 \cdot e_2 = e_4 \cdot e_5 = e_2$, $e_1 \cdot e_3 = e_4 \cdot e_6 = e_3$, we can do the re-factorization of obtained HNS by the subgroup $\{e_1, e_4\}$. As a result we will get another one HNS of the third dimension:

Table 7.

|  | $e_1$ | $e_2$ | $e_3$ |
|---|---|---|---|
| $e_2$ | | $\frac{1}{2}(e_1+e_3)$ | $\frac{1}{2}(e_1+e_2)$ |
| $e_3$ | | $\frac{1}{2}(e_1+e_2)$ | $\frac{1}{2}(e_1+e_3)$ |

On factorizing the hypercomplex system $\Gamma = N \cup \{0\}$ by the subgroup $A_6 = \{a_k : a_k = 6k, k \in Q\}$ with the convolution (11) we will get the finite hypercomplex system $G_6 = Q/A_6 = \{e_1, e_2, e_3, e_4, e_5, e_6\}$ of the sixth dimension. In this case, the multiplication of hypercomplex number system, taking into account (11) takes the form:

Table 8.

|  | $e_1$ | $e_2$ | $e_3$ | $e_4$ | $e_5$ | $e_6$ |
|---|---|---|---|---|---|---|
| $e_2$ | | $\frac{1}{2}(e_1+e_3)$ | $\frac{1}{2}(e_2+e_4)$ | $\frac{1}{2}(e_3+e_5)$ | $\frac{1}{2}(e_4+e_6)$ | $\frac{1}{2}(e_1+e_5)$ |
| $e_3$ | | $\frac{1}{2}(e_2+e_4)$ | $\frac{1}{2}(e_1+e_5)$ | $\frac{1}{2}(e_2+e_6)$ | $\frac{1}{2}(e_1+e_3)$ | $\frac{1}{2}(e_2+e_4)$ |
| $e_4$ | | $\frac{1}{2}(e_3+e_5)$ | $\frac{1}{2}(e_2+e_6)$ | $e_1$ | $e_2$ | $e_3$ |
| $e_5$ | | $\frac{1}{2}(e_4+e_6)$ | $\frac{1}{2}(e_1+e_3)$ | $e_2$ | $\frac{1}{2}(e_1+e_3)$ | $\frac{1}{2}(e_2+e_4)$ |
| $e_6$ | | $\frac{1}{2}(e_1+e_5)$ | $\frac{1}{2}(e_2+e_4)$ | $e_3$ | $\frac{1}{2}(e_2+e_4)$ | $\frac{1}{2}(e_1+e_5)$ |

Since there are such regularities as $e_1 \cdot e_1 = e_4 \cdot e_4 = e_1$, $e_1 \cdot e_2 = e_4 \cdot e_5 = e_2$, $e_1 \cdot e_3 = e_4 \cdot e_6 = e_3$ in current HNS's multiplication table, we can do re-factorization of $G_6$ by the subgroup $\{e_1, e_4\}$. As a result the HNS of the third dimension is got (Table 9).

Table 9.

|  | $e_1$ | $e_2$ | $e_3$ |
|---|---|---|---|
| $e_2$ | | $\frac{1}{2}(e_1+e_3)$ | $\frac{1}{2}(e_1+e_2)$ |
| $e_3$ | | $\frac{1}{2}(e_1+e_2)$ | $\frac{1}{2}(e_1+e_3)$ |

Likewise HNS of higher dimensions can be obtained. As we can see from the examples listed above, the results of factorization by the subgroups $A_6 = \{a_k : a_k = 6k, k \in Q\}, A_8 = \{a_k : a_k = 8k, k \in Q\}, A_{10} = \{a_k : a_k = 10k, k \in Q\}$ etc, that is, the subgroups, the elements of which are multiples of 2, can be re-factorized. Thus, we will get the HNS of two times less dimension.

## Conclusions

The method of conversion the infinite hypercomplex system to finite noncanonical hypercomplex number system by factorization methods is shown. Factorization of infinite hypercomplex number system by infinite subgroups is investigated with the following multiplication rules: $\sigma_n \cdot \sigma_m = \frac{1}{2}(\sigma_{n+m} + \sigma_{|n-m|})$. The obtained finite systems, starting with third dimension, are noncanonical. Obtained systems of even dimensions (2, 4, 6, 8…) can be re-factorized. As a result, the hypercomplex system of two times less dimension is got.